УДК 511, 517

# И. Ш. Джаббаров
## О показателе сходимости одного несобственного интеграла

**Аннотация**

В статье доказан новый результат о сходимости особого интеграла проблемы Терри. Результат основан на представление особого интеграла в виде некоторого поверхностного интеграла, распространенного на многообразие решений системы, связанного с рассматриваемой системой уравнений проблемы Терри. Получены новые оценки возникающих поверхностных интегралов, сводящие вопрос об оценке к исследованию некоторых операторов с дискретным спектром, с использованием «максимальных миноров».

**Ключевые слова:** поверхностные интегралы, тригонометрические интегралы, определитель Грама, алгебраические многообразия, неявные функции.

## Введение

В первые десятилетия XX столетия И. М. Виноградов разработал метод тригонометрических сумм, как средство решения широкого круга задач аналитической теории чисел. Предложенный им метод позволяет нетривиально оценить модули важных тригонометрических сумм, называемых суммами Вейля. Основой для оценок является теорема о среднем И. М. Виноградова ([1]).

Рассмотрим такую систему диофантовых уравнений:

$$x_1 + \cdots + x_k = y_1 + \cdots + y_k,$$
$$x_1^2 + \cdots + x_k^2 = y_1^2 + \cdots + y_k^2, \qquad (1)$$
$$\cdots \quad \cdots \quad \cdots$$
$$x_1^n + \cdots + x_k^n = y_1^n + \cdots + y_k^n,$$

где неизвестные $x_1, \cdots, x_k, y_1, \cdots, y_k$ принимают значения целых чисел от 1 до некоторого достаточно большого $P \geq 1$. Метод И. М. Виноградова основывается на оценке числа решений такой системы диофантовых уравнений. Проблема Терри представляет собой задачу об асимптотической формуле для числа решений системы (1) при $P \to \infty$. Относительно истории вопроса можно обратиться к [1-9]. В 1938 г. Хуа Ло-ген в работе [8] дал решение проблемы Терри на основе метода И. М. Виноградова. Пусть $I_{k,n}(P)$ обозначает число решений системы (1). Формула Хуа Ло-гена выглядит следующим образом:

$$I_{k,n}(P) = \sigma \theta_0 P^{2k-0,5(n^2+n)} + O\left(P^{2k-0,5(n^2+n)-\delta}\right),$$

где $\delta = \delta(n,k) > 0$ а $k$ имеет порядок $n^2 \ln n$, $\sigma$ - особый ряд (см. [9, стр.65]), $\theta_0$ - особый интеграл. Особый интеграл имеет вид (случай многочлена одного переменного):

$$\theta_0 = \theta_0(k) = \int\limits_{-\infty}^{\infty}\int\limits_{-\infty}^{\infty}\cdots\int\limits_{-\infty}^{\infty}\left|\int\limits_0^1 e^{2\pi i(\alpha_n x^n + \alpha_{n-1}x^{n-1}+\cdots+\alpha_1 x)}dx\right|^{2k} d\alpha_n d\alpha_{n-1}\cdots d\alpha_1.$$

Так, Хуа Ло-ген пришел к необходимости исследования показателя сходимости особого интеграла $\theta_0$ проблемы Терри.

***Определение.*** Число $\gamma > 0$ называется показателем сходимости несобственного интеграла $\theta_0(k)$, если он при вещественных $2k < \gamma$ расходится, а при $2k > \gamma$ сходится.

Хуа Ло-ген установил (см. [8]), что при $2k > 0,5n^2 + n$ особый интеграл сходится, указывая, при этом, что вопрос о точном значении показателя сходимости остается открытым. Последний вопрос был полностью решен авторами работы [4]. Одновременно, они выяснили, что значение показателя сходимости в случаях «полного» и «неполного» многочленов разные в следующем смысле. Многочлен $\alpha_n x^n + \alpha_{n-1}x^{n-1}+\cdots+\alpha_0 x$ они назвали неполным, если хотя бы один из коэффициентов $\alpha_i$ тождественно равен нулю (т. е. соответствующая степень переменной отсутствует в записи многочлена) и полным, в противном случае. Оказалось, что $\gamma = r + \cdots + m + n$ для неполных многочленов вида

$$\alpha_n x^n + \alpha_m x^m + \cdots + \alpha_r x^r; 1 \le r < \ldots < m < n$$

и $\gamma = 1 + 0,5(n^2 + n)$ для полных многочленов (заметим, что в полном случае

$$r + \cdots + m + n = 0,5(n^2 + n) < \gamma$$

и показатель сходимости на 1 больше суммы степеней одночленов).

С 1968 г. Архипов Г. И., Карацуба А. А. и Чубариков В. Н. начали создавать теорию кратных тригонометрических сумм, подобную теории Виноградова И. М.. Они обобщили теорему о среднем Виноградова И. М. на многомерный случай. В частности, авторы работ [3, 9] пришли к многомерному аналогу проблемы Терри и поставили вопрос о нахождении показателя сходимости особого интеграла многомерной проблемы.

Пусть

$$F(x,y) = \sum_{i=0}^{n}\sum_{j=0, i+j>0}^{m} \alpha_{ij} x^i y^j \quad (2)$$

некоторый многочлен от двух переменных $x$ и $y$ с действительными коэффициентами. Под особым интегралом двумерной проблемы Терри понимается интеграл

$$\theta_k = \int_{-\infty}^{\infty}\int_{-\infty}^{\infty}\cdots\int_{-\infty}^{\infty}\left|\int_0^1\int_0^1 e^{2\pi i F(x,y)}dxdy\right|^{2k} d\alpha_{10}d\alpha_{01}\cdots d\alpha_{mn}. \quad (3)$$

Показатель сходимости $\theta_k$ определяется как выше.

В общем многомерном случае существенное продвижение в рассматриваемых вопросах требует более глубокое изучение кратных тригонометрических интегралов и связанных с ними оценок мер многомерных областей ([10, 11]).

Вопрос об оценке сверху показателя сходимости рассматривался в работах ([4, 5, 7, 8, 9, 10]). В [9, стр. 50] доказана общая теорема, из которой для нашего двумерного случая следует оценка сверху $\gamma \leq N\max(n,m)$ где $N = (n+1)(m+1)-1$. В [10] приведена общая оценка сверху, которая для $\gamma$ дает следующую границу

$$\gamma \leq 2 + (n+m)(n+1)(m+1)/2$$

Эта оценка дает некоторое улучшение, когда разность $n-m$ велика по абсолютной величине. Рассматривался также и вопрос о точном значении показателя сходимости. В работах [13, 14] получены некоторые результаты в специальных случаях, очень близких к одномерному.

В настоящей работе мы, опираясь на результаты работ [10, 11, 15, 16] рассматриваем один модельный частный случай, когда сведение к одномерному случаю представляется невозможным. Мы ограничиваемся только верхней оценкой показателя сходимости.

Рассмотрим многочлен

$$F(x,y) = \alpha_1 x^3 + \alpha_2 x^2 y + \alpha_3 xy^2 + \alpha_4 y^3 + \alpha_5 x^2 + \alpha_6 xy + \alpha_7 y^2 + \alpha_8 x + \alpha_9 y, \quad (4)$$

состоящий из суммы трех однородных многочленов первой, второй и третьей степеней. Цель настоящей статьи состоит в доказательстве следующей теоремы.

**ТЕОРЕМА**. Интеграл

$$\int_{-\infty}^{\infty}\cdots\int_{-\infty}^{\infty}\left|\int_0^1\int_0^1 e^{2i\pi F(x,y)}dxdy\right|^{12} d\alpha_1\cdots d\alpha_9$$

сходится; здесь $F(x,y)$ - многочлен, определенный равенством (4).

Для сравнения заметим, что из результата кн. [9, стр. 51] сходимость получается только при замене показателя 12 на 24. А из результата нижней оценки работы [15] следует расходимость соответствующего интеграла, если 12 заменить числом 10. Ожидаемое точное значение показателя сходимости равно 11.

## 1. Вспомогательные леммы.

Для доказательства теоремы мы будем использовать некоторые вспомогательные утверждения. Эти утверждения содержатся в ниже сформулированных леммах.

*Лемма 1.* Пусть в ограниченной замкнутой жордановой области $\Omega$ $n$-мерного пространства $R^n$ задана непрерывная функция $f(\bar{x}) = f(x_1,...,x_n)$ и непрерывно-дифференцируемые функции $f_j(\bar{x}) = f_j(x_1,...,x_n)$, где $j = 1,...,r$, $r < n$, такие, что матрица Якоби

$$\frac{\partial(f_1,...,f_n)}{\partial(x_1,...,x_n)}$$

имеет всюду в $\Omega$ максимальный ранг. Пусть, далее $\bar{\xi}_0 = (\xi_1^0,...,\xi_r^0)$ - внутренняя точка образа отображения $\bar{x} \mapsto (f_1,...,f_r)$ и $\bar{x}_0$ - точка области $\Omega$ такая, что

$$f_1(\bar{x}_0) = \xi_1^0,..., f_0(\bar{x}_0) = \xi_r^0.$$

Тогда, всюду в некоторой окрестности точки $\bar{\xi}_0$ имеет место равенство

$$\frac{\partial^r}{\partial \xi_1 \cdots \partial \xi_r} \int_{\Omega(\bar{\xi})} f(\bar{x})d\bar{x} = \int_{M(\bar{\xi})} f(\bar{x}) \frac{ds}{\sqrt{G}},$$

где $\Omega(\bar{\xi})$ - подобласть в $\Omega$, определяемая системой неравенств $f_j(\bar{x}) \leq \xi_j$, $M(\bar{\xi})$ - поверхность, определяемая системой уравнений $f_j(\bar{x}) \leq \xi_j$ ($j = 1,...,r$), а $G$ - определитель Грама градиентов функций $f_j(\bar{x})$, т. е. $G = |(\nabla f_i, \nabla f_i)|$; на правой части стоит поверхностный интеграл первого рода (см. например [17, стр. 615], [18, стр. 263]).

*Следствие.* Пусть выполнены условия леммы 1. Тогда имеет место равенство

$$\int_{\Omega} f(\bar{x})d\bar{x} = \int_{m_1}^{M_1} \cdots \int_{m_r}^{M_r} du_1 \cdots du_r \int_M f(\bar{x}) \frac{ds}{\sqrt{G}},$$

где $m_j$ и $M_j$ - соответственно, минимальное и максимальное значения $f_j(\bar{x})$, $m_j \leq f_j(\bar{x}) \leq M_j$, $M = M(\bar{u})$ - поверхность в $\Omega$, определяемая системой уравнений $f_j = u_j, j = 1,...,r$, а $G$ - определитель Грама градиентов функций, определяющих $M$.

*Лемма 2.* Пусть, при условиях леммы 1, $\xi_1^0 = ... = \xi_r^0 = 0$ и поверхность $M$ определена системой уравнений

$$f_1(x_1,...,x_n) = 0,$$
$$\ldots$$
$$f_r(x_1,...,x_n) = 0,$$

причем, функции $f_j(\bar{x})$ непрерывно дифференцируемы в некоторой области $\Omega_0$, включающей внутри себя $\Omega$. Пусть $G = G(\bar{x})$ –определитель Грама градиентов функций $f_j(\bar{x})$, не равный нулю в $\Omega$. Пусть, далее, преобразование координат $\bar{x} = \bar{x}(\bar{\xi})$ взаимно – однозначно отображает некоторую область $\Omega'$ в $\Omega$, с неособой матрицей Якоби

$$Q = Q(\bar{\xi}) = \left(\frac{\partial x_i}{\partial \xi_j}\right)_{1 \leq i,j \leq n},$$

имеющей непрерывные в $\Omega$ элементы. Тогда, для любой непрерывной в $\Omega$ функции $f(\bar{x})$, имеет место формула

$$\int_M f(\bar{x})\frac{ds}{\sqrt{G}} = \int_{M'} |\det Q| f(\bar{x}(\bar{\xi}))\frac{d\sigma}{\sqrt{G'}}, G' = \det(JQ^tQ^tJ),$$

где $M'$-прообраз поверхности $M$ при этом преобразовании, $d\sigma$ - элемент площади в координатах $\bar{\xi}$, $J$ –матрица Якоби системы функций $f_j(\bar{x})$:

$$J = \frac{\partial(f_1,...,f_r)}{\partial(x_1,...,x_n)}$$

Доказательство этих утверждений даны в [12].

Нам понадобится многомерный аналог теоремы Арцела (см. [24, п. 546]):

*Лемма 3.* Пусть дана последовательность функций

$$f_n(\bar{x})(n = 1,2,...),$$

интегрируемых в произведении

$$K = [a_1,b_1] \times [a_2,b_2] \times \cdots \times [a_s,b_s]$$

и ограниченных в их совокупности

$$|f_n(\bar{x})| \leq L \; (\bar{x} \in K, n = 1,2,...).$$

Пусть для всех $\bar{x} \in K$ существует предел

$$\varphi(\bar{x}) = \lim_{n \to \infty} f_n(\bar{x}).$$

Если при любом $r, 0 \leq r \leq s-1$, функции $f_n(\bar{x})$ и функция $\varphi(\bar{x})$ интегрируемы на

$$[a_{r+1},b_{r+1}] \times \cdots \times [a_s,b_s],$$

то

$$\lim_{n \to \infty}\int_K f_n(\bar{x})d\bar{x} = \int_K \varphi(\bar{x})d\bar{x}.$$

*Доказательство.* Доказательство этой леммы проведем индукцией по $s$. При $s = 1$ получается одномерный случай теоремы Арцела (см. [24, п.546]). Пусть утверждение леммы справедливо при случае $s-1$.

При каждом $a_1 \le x_1 \le b_1$ полагаем

$$u_n(x_1) = \int_{a_2}^{b_2}\cdots\int_{a_s}^{b_s} f_n(\bar{x})dx_2\cdots dx_s.$$

По теореме о повторном интегрировании (см. [17, стр. 564 или стр. 568]) из условий леммы следует, что функция $u_n(x_1)$ интегрируема на $[a_1, b_1]$. Кроме того,

$$|u_n(x_1)| \le L\prod_{i=2}^{s}(b_s - a_s).$$

По индуктивному предположению под интегралом в выражении $u_n(x_1)$ можно переходить к пределу и полученная функция $\phi(x_1) = \int_{a_2}^{b_2}\cdots\int_{a_s}^{b_s}\varphi(\bar{x})dx_2\cdots dx_s$ интегрируема. Тогда, условия теоремы Арцела выполнены и

$$\lim_{n\to\infty}\int_{a_1}^{b_1} u_n(x_1)dx_1 = \int_{a_1}^{b_1}\lim_{n\to\infty} u_n(x_1)dx_1 = \int_{a_1}^{b_1}\phi(x_1)dx_1.$$

Равенство

$$\int_{a_1}^{b_1}\phi(x_1)dx_1 = \int_{K}\varphi(\bar{x})d\bar{x}$$

теперь следует из теоремы о повторном интегрировании.

*Лемма 4.* Пусть многочлен $F(x,y)$ определен равенством (4). Тогда, справедлива формула

$$\theta_k = (2\pi)^N \int_{\Pi} \frac{ds}{\sqrt{G_0}},$$

где поверхностный интеграл берется по поверхности $\Pi$, определяемой системой уравнений

$$\left.\begin{array}{r}x_1 + x_2 + x_3 + x_4 + x_5 + x_6 - x_7 - x_8 - x_9 - x_{10} - x_{11} - x_{12} = 0,\\ y_1 + y_2 + y_3 + y_4 + y_5 + y_6 - y_7 - y_8 - y_9 - y_{10} - y_{11} - y_{12} = 0,\\ x_1^2 + x_2^2 + x_3^2 + x_4^2 + x_5^2 + x_6^2 - x_7^2 - x_8^2 - x_9^2 - x_{10}^2 - x_{11}^2 - x_{12}^2 = 0,\\ x_1 y_1 + x_2 y_2 + x_3 y_3 + x_4 y_4 + x_5 y_5 + x_6 y_6 - x_7 y_7 - x_8 y_8 - x_9 y_9 - x_{10} y_{10} - x_{11} y_{11} - x_{12} y_{12} = 0,\\ y_1^2 + y_2^2 + y_3^2 + y_4^2 + y_5^2 + y_6^2 - y_7^2 - y_8^2 - y_9^2 - y_{10}^2 - y_{11}^2 - y_{12}^2 = 0,\\ x_1^3 + x_2^3 + x_3^3 + x_4^3 + x_5^3 + x_6^3 - x_7^3 - x_8^3 - x_9^3 - x_{10}^3 - x_{11}^3 - x_{12}^3 = 0,\\ x_1^2 y_1 + x_2^2 y_2 + x_3^2 y_3 + x_4^2 y_4 + x_5^2 y_5 + x_6^2 y_6 - x_7^2 y_7 - x_8^2 y_8 - x_9^2 y_9 - x_{10}^2 y_{10} - x_{11}^2 y_{11} - x_{12}^2 y_{12} = 0,\\ x_1 y_1^2 + x_2 y_2^2 + x_3 y_3^2 + x_4 y_4^2 + x_5 y_5^2 + x_6 y_6^2 - x_7 y_7^2 - x_8 y_8^2 - x_9 y_9^2 - x_{10} y_{10}^2 - x_{11} y_{11}^2 - x_{12} y_{12}^2 = 0,\\ y_1^3 + y_2^3 + y_3^3 + y_4^3 + y_5^3 + y_6^3 - y_7^3 - y_8^3 - y_9^3 - y_{10}^3 - y_{11}^3 - y_{12}^3 = 0\end{array}\right\} \quad (5)$$

в 24- мерном единичном кубе, а $G_0$ - определитель Грама градиентов, функций, стоящих на левых частях системы (5), т. е. $G_0 = \det(A_0 \cdot {}^t A_0)$ и

$$A_0 = \begin{pmatrix} 1 & 0 & 1 & 0 & \cdots & -1 & 0 & -1 & 0 \\ 0 & 1 & 0 & 1 & \cdots & 0 & -1 & 0 & -1 \\ 2x_1 & 0 & 2x_2 & 0 & \cdots & -2x_{11} & 0 & -2x_{12} & 0 \\ y_1 & x_1 & y_2 & x_2 & \cdots & -y_{11} & -x_{11} & -y_{12} & -x_{12} \\ 0 & 2y_1 & 0 & 2y_2 & \cdots & 0 & -2y_{11} & 0 & -2y_{12} \\ 3x_1^2 & 0 & 3x_2^2 & 0 & \cdots & -3x_{11}^2 & 0 & -3x_{12}^2 & 0 \\ 2x_1 y_1 & x_1^2 & 2x_2 y_2 & x_2^2 & \cdots & -2x_{11} y_{11} & -x_{11}^2 & -2x_{12} y_{12} & -x_{12}^2 \\ y_1^2 & 2x_1 y_1 & y_2^2 & 2x_2 y_2 & \cdots & -y_{11}^2 & -2x_{11} y_{11} & -y_{12}^2 & -2x_{12} y_{12} \\ 0 & 3y_1^2 & 0 & 3y_2^2 & \cdots & 0 & -3y_{11}^2 & 0 & -3y_{12}^2 \end{pmatrix};$$

при этом, поверхностный интеграл, в силу того, что подынтегральное выражение не определено для тех значений $\bar{x} \in [0,1]^{24}$, где $G_0 = 0$, определяется в несобственном смысле, как предел

$$\lim_{\eta \to 0} \int_{M_\eta \times M_\eta} \frac{ds}{\sqrt{G_0}},$$

где множество $M_\eta$ определено ниже (см. стр. 8).

Прежде, чем доказать лемму 4, заметим, что несобственный поверхностный интеграл определяется переходом к пределу по симметрическому прямому произведению, подобно главному значению интеграла по Коши в одномерном случае. Если вместо симметрического прямого произведения (которое определяется рассмотрением не системы функций из (7), а системы функций, стоящих на левых частях системы (5)) взять произвольные области, определяемые неравенствами вида $G_0 \geq \eta$, то утверждение теоремы может не выполнятся, из-за того, что в таком случае сходимость поверхностного интеграла, как становится ясным ниже из доказательства леммы 1, не равносильно сходимости особого интеграла.

*Доказательство леммы 4.* Возведем внутренний тригонометрический интеграл в равенстве теоремы в 6-ю степень:

$$\left( \int_0^1 \int_0^1 e^{2i\pi F(x,y)} dx dy \right)^6 = \int_0^1 \int_0^1 \cdots \int_0^1 \int_0^1 e^{2i\pi (F(x_1,y_1) + \cdots + F(x_6,y_6))} dx_1 dy_1 \cdots dx_6 dy_6.$$

Подставляя (4) в эту формулу мы получим:

$$\left( \int_0^1 \int_0^1 e^{2i\pi F(x,y)} dx dy \right)^6 = \int_0^1 \int_0^1 \cdots \int_0^1 \int_0^1 e^{2i\pi (\alpha_1 u_1 + \cdots + \alpha_9 u_9)} dx_1 dy_1 \cdots dx_6 dy_6, \quad (6)$$

где для краткости мы ввели обозначения:

$$u_j = f_j = f_j(x_1, y_1, \ldots, x_6, y_6); j = 1, \ldots, 9,$$

при этом,

$$f_1 = x_1 + x_2 + x_3 + x_4 + x_5 + x_6,$$

$$f_2 = y_1 + y_2 + y_3 + y_4 + y_5 + y_6,$$

$$f_3 = x_1^2 + x_2^2 + x_3^2 + x_4^2 + x_5^2 + x_6^2,$$

$$f_4 = x_1 y_1 + x_2 y_2 + x_3 y_3 + x_4 y_4 + x_5 y_5 + x_6 y_6,$$

$$f_5 = y_1^2 + y_2^2 + y_3^2 + y_4^2 + y_5^2 + y_6^2, \qquad (7)$$

$$f_6 = x_1^3 + x_2^3 + x_3^3 + x_4^3 + x_5^3 + x_6^3,$$

$$f_7 = x_1^2 y_1 + x_2^2 y_2 + x_3^2 y_3 + x_4^2 y_4 + x_5^2 y_5 + x_6^2 y_6,$$

$$f_8 = x_1 y_1^2 + x_2 y_2^2 + x_3 y_3^2 + x_4 y_4^2 + x_5 y_5^2 + x_6 y_6^2,$$

$$f_9 = y_1^3 + y_2^3 + y_3^3 + y_4^3 + y_5^3 + y_6^3.$$

Сначала заметим, что множество особых точек отображения $f:[0,1]^{12} \to [0,6]^9$, является замкнутым множеством в $I = [0,1]^{12}$ и потому имеет нулевую Жорданову меру, если доказать, что $G = 0$ на подмногообразии куба $I$ лебеговой меры нуль. Это будет установлено ниже. Пусть $M = \{(x_1, y_1,..., x_6, y_6) \in [0,1]^{12} \mid G = 0\}$. Взяв положительное число $\varepsilon > 0$ мы можем покрыть множество $M$ конечным семейством открытых кубов, суммарной меры $\varepsilon > 0$. Обозначим $E$ объединение кубов этого семейства. Пусть $\eta > 0$ минимальное значение функции $G = G(\bar{x})$ в $I \setminus E$. Имеем:

$$\int_0^1 \int_0^1 \cdots \int_0^1 \int_0^1 e^{2i\pi(\alpha_1 u_1 + \cdots + \alpha_9 u_9)} dx_1 dy_1 \cdots dx_6 dy_6 = \int_{I \setminus E} + \int_E .$$

Очевидно, второй интеграл по абсолютной величине не превосходит $\varepsilon$. Следовательно, при $\varepsilon \to 0$ первый интеграл на правой части стремится к значению интеграла на левой части. Тогда, если ввести в рассмотрение множества

$$M_\eta = \{(x_1, y_1,..., x_6, y_6) \in [0,1]^{12} \mid G \geq \eta > 0\},$$

то получим

$$\int_{I \setminus E} e^{2i\pi(\alpha_1 u_1 + \cdots + \alpha_9 u_9)} dx_1 dy_1 \cdots dx_6 dy_6 = \int_{M_\eta} e^{2i\pi(\alpha_1 u_1 + \cdots + \alpha_9 u_9)} dx_1 dy_1 \cdots dx_6 dy_6 + \delta; |\delta| \leq \varepsilon .$$

Теперь применяя следствие к лемме 1, под интегралом на правой части произведем замену переменных $u_j = f_j; j = 1,...,9$. Тогда получим:

$$\int_0^1 \int_0^1 \cdots \int_0^1 \int_0^1 e^{2i\pi(\alpha_1 u_1 + \cdots + \alpha_9 u_9)} dx_1 dy_1 \cdots dx_6 dy_6 =$$

$$= \int_0^6 \cdots \int_0^6 \left( \lim_{\eta \to 0} \int_{\Pi(\bar{u}) \cap M_\eta} \frac{ds}{\sqrt{G}} \right) e^{2i\pi(\alpha_1 u_1 + \cdots + \alpha_9 u_9)} du_1 \cdots du_9 ,$$

где $\Pi(\bar{u})$, при $\bar{u} = (u_1,...,u_9) \in [0,6]^9$, обозначает поверхность в $[0,1]^{12}$, определяемую системой уравнений $f_j = u_j; j = 1,...,9$.

*Замечание.* На правой части последнего равенства стоит несобственный интеграл, который понимается (по определению) как предел

$$\lim_{\eta \to 0} \int_0^6 \cdots \int_0^6 \left( \int_{\Pi(\bar{u}) \cap M_\eta} \frac{ds}{\sqrt{G}} \right) e^{2i\pi(\alpha_1 u_1 + \cdots + \alpha_9 u_9)} du_1 \cdots du_9 .$$

Поэтому, выражение $\lim_{\eta \to 0} \int_{\Pi(\bar{u}) \cap M_\eta} \frac{ds}{\sqrt{G}}$ мы рассматриваем, как функцию от $\bar{u} = (u_1,...,u_9)$, которая в силу положительности и монотонности имеет вполне определенное, конечное или бесконечное значение.

Левая часть полученного выше равенства представляет собой функцию от $\bar{\alpha} \in R^9$, а правая часть - преобразование Фурье поверхностного интеграла, зависящего от $\bar{u} = (u_1,...,u_9) \in [0,6]^9$. По теореме Планшереля принадлежность одной из них Лебеговому классу $L_2$ влечет принадлежность этому же классу и второй функции, при этом

$$\int_{-\infty}^\infty \cdots \int_{-\infty}^\infty \left| \int_0^1 \int_0^1 e^{2i\pi F(x,y)} dx dy \right|^{12} d\alpha_1 \cdots d\alpha_9 = \int_0^6 \cdots \int_0^6 \left( \lim_{\eta \to 0} \int_{\Pi(\bar{u}) \cap M_\eta} \frac{ds}{\sqrt{G}} \right)^2 du_1 \cdots du_9 .$$

Поэтому, для сходимости интеграла на левой части необходимо и достаточно сходимость интеграла на правой части. Докажем, что последнее равенство записывается в виде

$$\int_{-\infty}^\infty \cdots \int_{-\infty}^\infty \left| \int_0^1 \int_0^1 e^{2i\pi F(x,y)} dx dy \right|^{12} d\alpha_1 \cdots d\alpha_7 = \lim_{\eta \to 0} \int_0^6 \cdots \int_0^6 \left( \int_{\Pi(\bar{u}) \cap M_\eta} \frac{ds}{\sqrt{G}} \right)^2 du_1 \cdots du_9 . \qquad (8)$$

Пусть $D$ - некоторая замкнутая подобласть единичного куба, где $G > 0$. Тогда, по лемме 1 имеем:

$$\int_{\bar{x} \in D, \Pi(\bar{u})} \frac{ds}{\sqrt{G}} = \lim_{h \to 0} \frac{1}{(2h)^9} \int_{u_j - h < f_j < u_j + h} d\bar{x} ; \quad \bar{x} = (x_1,...,y_6). \qquad (9)$$

Разобьем область $M_\eta$, как при доказательстве леммы 1, п. 1, гл. 1, на замкнутые части $D'_\nu, \nu = 1,...,l$, в каждой из которых некоторый минор матрицы $A_0$ принимает, по модулю, максимальные значения. Тогда, и поверхностный интеграл разбивается на сумму интегралов, у которых подынтегральная функция ограничена. Обозначим левую часть (9) $\varphi_\eta(\bar{u})$. Докажем, что она ограничена числом не зависящим от $h$. Учитывая проведенное разбиение и обозначая модуль максимального минора как $|J|$, получим $\eta \leq G \leq C_{12}^9 |J|^2 = 220|J|^2$. Переходя к кратному интегралу, используя выражение для элемента площади, можем писать:

$$\int_{\bar{x}\in D, \Pi(\bar{u})} \frac{ds}{\sqrt{G}} \leq \frac{25}{\sqrt{\eta}} \int_0^1 \cdots \int_0^1 d\xi_1 \cdots d\xi_5 \leq \frac{25}{\sqrt{\eta}}.$$

Тогда, условия леммы 3 выполнены. Далее, по лемме 1 и следствию к ней

$$\int_0^k \cdots \int_0^k \varphi_D(\bar{u})\varphi_D(\bar{u})d\bar{u} = \int_0^k \cdots \int_0^k \varphi_D(\bar{u})\left(\lim_{h\to 0}\frac{1}{(2h)^9}\sum_\nu \int_{u_j-h<f_j<u_j+h, \bar{x}\in D_\nu} d\bar{x}\right)d\bar{u}.$$

Придадим $h$ произвольные значения $h = h_n$ такие, что $h_n \to 0$ и применим лемму 3 к интегралу на правой части последнего равенства. Тогда получим:

$$\int_0^k \cdots \int_0^k \varphi_D(\bar{u})\left(\lim_{n\to\infty}\frac{1}{(2h_n)^9} \int_{u_j-h_n<f_j<u_j+h_n, \bar{x}\in D_\nu} d\bar{x}\right)d\bar{u} =$$

$$= \lim_{n\to\infty}\frac{1}{(2h_n)^9}\int_0^k \cdots \int_0^k \varphi_D(\bar{u})\left(\int_{u_j-h_n<f_j<u_j+h_n, \bar{x}\in D_\nu} d\bar{x}\right)d\bar{u}.$$

Суммируя по $\nu$, будем иметь:

$$\int_0^k \cdots \int_0^k \varphi_D(\bar{u})\varphi_D(\bar{u})d\bar{u} = \lim_{h\to 0}\frac{1}{(2h)^9}\int_0^k \cdots \int_0^k \int_{\bar{x}'\in D,\Pi'(\bar{u})}\left(\int_{u_j-h<f_j<u_j+h,} d\bar{x}\right)\frac{ds'}{\sqrt{G'}}d\bar{u}, \quad (10)$$

здесь $ds'$ означает элемент площади поверхности $\Pi'(\bar{u})$, определяемой в $D$ системой уравнений $f_j(\bar{x}') = u_j; j=1,...,9$ причем, $G'$ имеет аналогичный смысл. Рассмотрим, при фиксированном $h$, внутренний интеграл в последней цепочке равенств (10), т. е. интеграл

$$\int_{\bar{x}'\in D,\Pi'(\bar{u})}\left(\int_{u_j-h<f_j<u_j+h} d\bar{x}\right)\frac{ds'}{\sqrt{G'}}.$$

Для точек $\bar{x}' \in [0,1]^{12}$ мы определяем функцию $f(\bar{x}')$, полагая ее значение, для каждого $\bar{x}' \in \Pi'(\bar{u})$, равным внутреннему интегралу. Докажем, что функция, $f(\bar{x}')$, определенная таким образом, непрерывна на $D$. Пусть $\bar{x}_1', \bar{x}_2' \in D$, $\bar{x}_1' = (x_{11}, y_{11},...,x_{16}, y_{16})$, $\bar{x}_{21}' = (x_{21}, y_{21},...,x_{26}, y_{26})$ $\sum[(x_{1i}'-x_{2i}')^2 + (y_{1i}'-y_{2i}')^2] \leq \varepsilon$, при заданном $\varepsilon > 0$. Тогда, полагая $u_j^1 = f_j(\bar{x}_1')$, $u_j^2 = f_j(\bar{x}_2')$ (здесь мы пользуемся верхним индексированием), согласно формуле о конечных приращениях будем иметь при некоторых $\bar{\theta}$ и $\bar{\eta}$:

$$|u_j^1 - u_j^2| = \left|\sum_{r=1}^6 \left(\frac{\partial f_j}{\partial x_r}(\bar{x}_1' + \bar{\theta})(x_{1r}' - x_{2r}') + \frac{\partial f_j}{\partial x_r}(\bar{y}_1' + \bar{\eta})(y_{1r}' - y_{2r}')\right)\right| \leq 3\sqrt{12\varepsilon},$$

если $\sum[(x'_{1i} - x'_{2i})^2 + (y'_{1i} - y'_{2i})^2] \leq \varepsilon$. Следовательно, находим:

$$\left| f(\bar{x}'_1) - f(\bar{x}'_2) \right| = \left| \int\limits_{u^1_j - h < f_j < u^1_j + h} d\bar{x} - \int\limits_{u^2_j - h < f_j < u^2_j + h} d\bar{x} \right|.$$

Выражение на правой части представляет собой абсолютное значение разности объемов прообразов двух кубов, с достаточно близкими центрами. Из геометрических соображений ясно, что оно не превосходит объемов прообразов прямоугольных параллелепипедов, примыкающих к граням кубов и окаймляющих их общую часть. Число этих параллелепипедов равно числу граней одного куба, т.е. $2^{12}$. Таким образом,

$$\left| f(\bar{x}'_1) - f(\bar{x}'_2) \right| \leq 2^{12} \max_j \left| \int\limits_{u^1_j - h - 6\sqrt{3\varepsilon} < f_j < u^1_j - h + 6\sqrt{3\varepsilon}} d\bar{x} + \int\limits_{u^2_j - h - 6\sqrt{3\varepsilon} < f_j < u^2_j + h + 6\sqrt{3\varepsilon}} d\bar{x} \right|.$$

Эти интегралы оцениваются одинаково. Поскольку, независимые переменные меняются в пределах [0,1], используя выражение для элемента площади, найденное при доказательстве леммы 1, находим:

$$\int\limits_{u^1_j - h - 6\sqrt{3\varepsilon} < f_j < u^1_j - h + 6\sqrt{3\varepsilon}} d\bar{x} \leq \int\limits_{u^1_j - h - 6\sqrt{3\varepsilon}}^{u^1_j - h + 6\sqrt{3\varepsilon}} du_1 \int\limits_0^6 \cdots \int\limits_0^6 du_2 \cdots du_9 \int\limits_{\Pi(\bar{u})} \frac{ds}{\sqrt{G}} \leq 25 \cdot 6^9 \sqrt{\frac{3\varepsilon}{\eta}}.$$

Учитывая предыдущее неравенство, в силу произвольности $\varepsilon$, получаем непрерывность функции $f(\bar{x}')$. Таким образом,

$$\int\limits_{\Pi'(\bar{u})} \left( \int\limits_{u_j - h < f_j < u_j + h} d\bar{x} \right) \frac{ds'}{\sqrt{G'}} = \int\limits_{\Pi'(\bar{u})} f(\bar{x}') \frac{ds'}{\sqrt{G'}}.$$

Применяя следствие к лемме 1, находим:

$$\int\limits_0^6 \cdots \int\limits_0^6 \int\limits_{\bar{x}' \in D, \Pi'(\bar{u})} f(\bar{x}') \frac{ds'}{\sqrt{G'}} d\bar{u} = \int\limits_0^1 \cdots \int\limits_0^1 f(\bar{x}') d\bar{x}' =$$

$$= \int\limits_0^1 dx'_1 \cdots \int\limits_0^1 dy'_6 \int\limits_{u_j - h \leq f_j \leq u_j - h, j=1,\ldots 9} dx_1 \cdots dy_6.$$

Учитывая, что $f_j(\bar{x}') = u_j; j = 1,\ldots,9$, при новых обозначениях $x_{6+s} = x'_s, y_{6+s} = y'_s$, из (10) окончательно получаем:

$$\int\limits_0^k \cdots \int\limits_0^k \varphi_D(\bar{u}) \varphi_D(\bar{u}) d\bar{u} =$$

$$= \lim_{h \to 0} \frac{1}{(2h)^9} \int\limits_0^1 dx'_1 \cdots \int\limits_0^1 dy'_6 \int\limits_{u_j - h \leq f_j \leq u_j - h, j=1,\ldots 9} dx_1 \cdots dy_6 = \int\limits_{\bar{x} \in D \times D, \Pi_0} \frac{ds}{\sqrt{G}}; \bar{x} = (x_1, \ldots, y_{12}).$$

где $П_0$ является многообразием решений системы (5), а $G_0$ - определителем Грама градиентов функций $f_j - f_j'$, стоящих на левых частях системы (5). Из соотношения для определителя Грама (см. [22, стр. 245]) следует, что $G_0 \geq G + G' > 0$ в $D \times D$. Для завершения доказательства остается переходить к пределу при $\eta \to 0$. Лемма 4 доказана.

*Лемма 5.* Пусть в пространстве $R^n$ задана полиномиальная система уравнений
$$f_1(\bar{x}) = 0,..., f_k(\bar{x}) = 0; k < n,$$
с многочленами на левых частях, никакой из которых не делится на другой. Тогда, алгебраическое множество, определяемое этой системой уравнений, является подмножеством алгебраического множества, состоящего из объединения конечного числа (зависящего от числа уравнений и их степеней) поверхностей размерности не выше $n-k$.

*Доказательство.* Доказательство леммы 5 проведем методом математической индукции относительно количества уравнений системы. При $k=1$ ее утверждение следует из теоремы 3.2, [25]. Предположим, что утверждение леммы справедливо для систем с $k-1$ уравнений. Докажем ее для систем с $k$ уравнений. Возьмем первое уравнение системы. Точки, являющиеся решениями системы служат, также общим решением всех систем уравнений
$$f_1(\bar{x}) = 0, f_j(\bar{x}) = 0; j = 2,...,n,$$
первым из которых является первое уравнение данной системы. Обозначим результаты таких систем через $K_j(\bar{x}')$, $\bar{x}' = (x_1,...,x_{n-1})$. Если данная система имеет решение $\bar{x}$, то точка $\bar{x}' = (x_1,...,x_{n-1})$ является решением системы
$$K_j(\bar{x}') = 0, j = 2,...,k;$$
при этом, $f_1(\bar{x}', x_n) = 0$. По теореме 3.2, [25], все решения последнего уравнения описываются при помощи конечного числа равенств вида $x_n = \omega(\bar{x}')$. По индуктивному предположению, все решения последней системы расположены на конечного числа замкнутых поверхностях вида
$$x_{n-j} = \psi_{n-j}(x_1,...,x_{n-k}), j = 1,...,k-1.$$
Таким образом, все решения данной системы расположены на конечного числа замкнутых поверхностях размерности не выше $n-k$, вида
$$x_n = \omega(x_1,...,x_{n-k}, \psi_{n-k+1}(x_1,...,x_{n-k}),...,\psi_{n-1}(x_1,...,x_{n-k})).$$
Лемма 5 доказана.

*Следствие.* Рассмотрим в пространстве $R^n$ произвольную систему полиномиальных неравенств

$$f_1(\bar{x}) \geq 0,$$
$$\vdots$$
$$f_n(\bar{x}) \geq 0,$$

в произвольном ограниченном параллелепипеде. Тогда множество решений этой системы, заключенное в рассматриваемом параллелепипеде состоит из объединения конечного числа замкнутых односвязных областей, причем число этих областей зависит лишь от степеней левых частей неравенств данной системы.

## 2. Доказательство теоремы

Для доказательства теоремы мы должны, показать, что отображение определенное системой равенств (7) является неособым всюду в кубе $[0,1]^{12}$, за исключением некоторого жорданового подмножества жордановой меры нуль. Для этого мы должны, согласно замечаниям при доказательстве леммы 4, установить, что $G \neq 0$ вне подмногообразии куба $[0,1]^{12}$ жордановой меры нуль. Далее, для завершения доказательства теоремы остается показать, что поверхностный интеграл леммы 4 сходится.

Докажем первое утверждение. Рассмотрим матрицу

$$A = \begin{pmatrix} 1 & 0 & 1 & 0 & \cdots & 1 & 0 & 1 & 0 \\ 0 & 1 & 0 & 1 & \cdots & 0 & 1 & 0 & 1 \\ 2x_1 & 0 & 2x_2 & 0 & \cdots & 2x_5 & 0 & 2x_6 & 0 \\ y_1 & x_1 & y_2 & x_2 & \cdots & y_5 & x_5 & y_6 & x_6 \\ 0 & 2y_1 & 0 & 2y_2 & \cdots & 0 & 2y_5 & 0 & 2y_6 \\ 3x_1^2 & 0 & 3x_2^2 & 0 & \cdots & 3x_5^2 & 0 & 3x_6^2 & 0 \\ 2x_1 y_1 & x_1^2 & 2x_2 y_2 & x_1^2 & \cdots & 2x_5 y_5 & x_5^2 & 2x_6 y_6 & x_6^2 \\ y_1^2 & 2x_1 y_1 & y_2^2 & 2x_2 y_2 & \cdots & y_5^2 & 2x_5 y_5 & y_6^2 & 2x_6 y_6 \\ 0 & 3y_1^2 & 0 & 3y_2^2 & \cdots & 0 & 3y_5^2 & 0 & 3y_6^2 \end{pmatrix},$$

и возьмем минор, образованный столбцами с номерами 1, 3, 5, 7, 9, 6, 8, 10, 12:

$$\begin{vmatrix} 1 & 1 & 1 & 1 & 1 & 0 & 0 & 0 & 0 \\ 0 & 0 & 0 & 0 & 0 & 1 & 1 & 1 & 1 \\ 2x_1 & 2x_2 & 2x_3 & 2x_4 & 2x_5 & 0 & 0 & 0 & 0 \\ y_1 & y_2 & y_3 & y_4 & y_5 & x_3 & x_4 & x_5 & x_6 \\ 0 & 0 & 0 & 0 & 0 & 2y_3 & 2y_4 & 2y_5 & 2y_6 \\ 3x_1^2 & 3x_2^2 & 3x_3^2 & 3x_4^2 & 3x_5^2 & 0 & 0 & 0 & 0 \\ 2x_1 y_1 & 2x_2 y_2 & 2x_3 y_3 & 2x_4 y_4 & 2x_5 y_5 & x_3^2 & x_4^2 & x_5^2 & x_6^2 \\ y_1^2 & y_2^2 & y_3^2 & y_4^2 & y_5^2 & 2x_3 y_3 & 2x_4 y_4 & 2x_5 y_5 & 2x_6 y_6 \\ 0 & 0 & 0 & 0 & 0 & 3y_3^2 & 3y_4^2 & 3y_5^2 & 3y_6^2 \end{vmatrix}.$$

Поменяя местами строк перепишем его в виде:

$$-\begin{vmatrix} 3x_1^2 & 3x_2^2 & 3x_3^2 & 3x_4^2 & 3x_5^2 & 0 & 0 & 0 & 0 \\ y_1^2 & y_2^2 & y_3^2 & y_4^2 & y_5^2 & 2x_3 y_3 & 2x_4 y_4 & 2x_5 y_5 & 2x_6 y_6 \\ 2x_1 & 2x_2 & 2x_3 & 2x_4 & 2x_5 & 0 & 0 & 0 & 0 \\ y_1 & y_2 & y_3 & y_4 & y_5 & x_3 & x_4 & x_5 & x_6 \\ 2x_1 y_1 & 2x_2 y_2 & 2x_3 y_3 & 2x_4 y_4 & 2x_5 y_5 & x_3^2 & x_4^2 & x_5^2 & x_6^2 \\ 1 & 1 & 1 & 1 & 1 & 0 & 0 & 0 & 0 \\ 0 & 0 & 0 & 0 & 0 & 2y_3 & 2y_4 & 2y_5 & 2y_6 \\ 0 & 0 & 0 & 0 & 0 & 1 & 1 & 1 & 1 \\ 0 & 0 & 0 & 0 & 0 & 3y_3^2 & 3y_4^2 & 3y_5^2 & 3y_6^2 \end{vmatrix}.$$

Обозначим через $B$ матрицу порядка 5, стоящую в левом верхнем углу. Легко заметить, что при раскрытии определителя этой матрицы мы получаем сумму 5! попарно различных одночленов, которые не могут исчезать приведением подобных членов. Поэтому, мы получим ненулевой многочлен 8-й степени от 10-ти переменных. Множество нулей этого многочлена является подмногообразием размерности 9 и имеет меру нуль в 10-ти мерном пространстве (оно является проекцией некоторого подмногообразия размерности 11 в $[0,1]^{12}$). Рассмотрим точки $[0,1]^{12}$ лежащие вне этого подмногообразия, где $|B| \neq 0$. Обозначим через $C$ - матрицу размера $4 \times 5$, расположенную, под матрицей $B$, с единственной ненулевой (первой) строкой, а через $D$ и $E$ - соответственно матриц размеров $5 \times 4$ и $4 \times 4$ так, чтобы определитель можно было представить в виде следующего блочного определителя:

$$\begin{vmatrix} B & D \\ C & E \end{vmatrix} = |B| \left| E - |B|^{-1} CB^{\vee} D \right|$$

(см. формулу [21, стр. 127, задача 46]). Произведение $CB^{\vee}$ имеет только одну ненулевую (первую) строку. Поэтому, то же самое верно для матрицы $CB^{\vee}D$. Тогда, первая строка матрицы

$$E - |B|^{-1} CB^{\vee} D$$

состоит (так как первая строка $E$ - нулевая) из первой строки матрицы $CB^{\vee}D$, помноженную на $-|B|^{-1}$, а оставшиеся строки такие же как у матрицы $E$. Вынося этот множитель за знак определителя, мы получим:

$$\begin{vmatrix} B & D \\ C & E \end{vmatrix} = \left| E - CB^{\vee} D \right|.$$

Обозначим суммы алгебраических дополнений элементов первой, второй и т. д. строк матрицы $B$, соответственно, $b_1, b_2, b_3, b_4, b_5$. Первая строка матрицы $CB^{\vee}D$ равна

$$\sum_{j=1}^{5} b_j D^j,$$

где $D^j$ есть $j$-я строка матрицы $D$. Итак,

$$\begin{vmatrix} B & D \\ C & E \end{vmatrix} = \begin{vmatrix} \varphi(x_3,y_3) & \varphi(x_4,y_4) & \varphi(x_5,y_5) & \varphi(x_6,y_6) \\ 2y_3 & 2y_4 & 2y_5 & 2y_6 \\ 1 & 1 & 1 & 1 \\ 3y_3^2 & 3y_4^2 & 3y_5^2 & 3y_6^2 \end{vmatrix},$$

причем, $\varphi(x_j,y_j) = 2x_j y_j b_2 + 2y_j b_4 + x_j^2 b_5$. Так как многочлены $b_1, b_2, b_3, b_4, b_5$ не содержат $x_6$, то последний определитель в своем раскрытии содержит $x_6^2$ только в произведении

$$\varphi(x_6, y_6) \begin{vmatrix} 2y_3 & 2y_4 & 2y_5 \\ 1 & 1 & 1 \\ 3y_3^2 & 3y_4^2 & 3y_5^2 \end{vmatrix},$$

т. е. в произведении $b_5 x_6^2 (y_5 - y_3)(y_4 - y_3)(y_4 - y_5)$. Следовательно, $G \neq 0$ вне подмногообразий, определенных равенствами $b_5 = 0$, $x_6 = 0$, $y_5 = y_3, y_4 = y_3, y_4 = y_5$. Как уже отмечено выше, подмногообразие $b_5 = 0$ имеет размерность 11 в $[0,1]^{12}$, то тем самым мы доказали, что равенство $G = 0$ может выполняться на подмножестве куба $[0,1]^{12}$ жордановой меры нуль.

Теперь покажем, что поверхностный интеграл леммы 4 сходится в смысле указанном в утверждении этой леммы.

1) Пусть $G = \max_{\bar{x} \in \Omega} G_0(\bar{x})$. Тогда, (см. [22, задача 13.30]) из следующего неравенства видно, что $G$ ограничено:

$$G_0 = \det(A_0 \cdot {}^t A_0) \leq G \leq 2^{12} \cdot 12^{10} \leq 2^{52}.$$

Положим $\Omega_p = \{\bar{x} \in \Omega \mid 2^{52-2p} \leq G_0 \leq 2^{54-2p}\}$ при $p = 1,2,\ldots$. Обозначим $\Pi_p$-часть поверхности $\Pi$, лежащую в $\Omega_p$. Если взять произвольную замкнутую область, лежащую в области, где $G > 0$, то в этой области $G_0$ достигнет свою нижнюю грань. Поэтому, такая область лежит в объединении конечного числа областей вида $\Omega_p$. Тогда, сходимость $\theta_k$ будет следовать из сходимости ряда $\sum_{p=1}^{\infty} I_p$, где $I_p = \int_{\Pi_p} \frac{ds}{\sqrt{G_0}}$.

2) Пусть $p$ фиксирован. Оценим интеграл $I_p$. Для этого произведем в интеграле $I_p$ замену переменных по формулам: $x_j = (2h)^\alpha u_j$, $y_j = (2h)^\alpha v_j$ $j = 1,\ldots,12$, $h = 2^{27-p}$, $\alpha = 1/11$. Матрица Якоби $Q$ этого преобразования диагональная матрица порядка 24, определитель которой равен $(2h)^{24/11}$. Каждая строка матрицы $A_0 = A_0(\bar{x})$ содержит одночлены одинаковой степени и поэтому, после подстановки $x_j = (2h)^{\alpha} u_j$, $y_j = (2h)^{\alpha} v_j$ матрица $A_0 = A_0(\bar{x})$ будет иметь общие множители в строках с номерами 3-9. Их можно вынести за знак определителя Грама, возводя предварительно в квадрат. Из вида матрицы $Q$ следует, что $A_0 Q = (2h)^{\alpha} A_0$ и поэтому,

$$G'_0 = \det(A_0 Q \cdot {}^t Q {}^t A_0) = (2h)^{18/11} \det(A_0 \cdot {}^t A_0).$$

По лемме 2, $G'_0 = \det(A_0 Q \cdot {}^t Q {}^t A_0)$; следовательно,

$$\int_{\Pi_p} \frac{ds}{\sqrt{G_0}} = \int_{\Pi'_p} |\det Q| \frac{d\sigma}{\sqrt{G'_0}}$$

где $\Pi'_p$ - прообраз поверхности $\Pi_p$, при произведенной замене переменных, $d\sigma$ - элемент площади на $\Pi'_p$. Имеем:

$$\det(A_0 \cdot {}^t A_0) = (2h)^2 G''_0,$$

причем $G''_0$ имеет такой же вид, что и $G_0$, но только переменные интегрирования $u_j, v_j$, меняются в пределах $0 \leq u_j, v_j \leq (2h)^{-\alpha}$. Тогда,

$$I_p = (2h)^{4/11} \int_{\Pi'_p} \frac{d\sigma}{\sqrt{G''_0}}; \quad 0,5 < \sqrt{G''_0} \leq 1, \quad (12)$$

При этом, матрицу $A_0(\bar{x})$ леммы 4, после сокращения всех элементов строк на общие множители, снова обозначаем как $A_0 = A_0(\bar{u})$, $\bar{u} = (u_1, v_1, ..., u_{12}, v_{12})$.

3) Обозначим $\Omega(h)$ ту часть куба $[0, (2h)^{-1/11}]^{24} = [0, 2^{p-28}]^{24}$, где выполняется неравенство $0,5 \leq \sqrt{G''_0} \leq 1$. Для того, чтобы оценить поверхностный интеграл на правой части (12), мы разобьем его, предварительно, на такие части, чтобы величину интеграла по каждой такой части можно было оценить, используя подходящую проекцию, сопоставимую по величине с куском соответствующей части поверхности. Эти части определяются максимальными минорами матрицы Якоби $A_0$ системы леммы 4. Разобьем $\Omega(h)$ на не более, чем $t = \binom{24}{9}$ подобластей $\Omega^{(\nu)}, \nu = 1, ..., t$, в каждой из которых один из миноров матрицы $A_0 = A_0(\bar{u})$, имеет, по модулю, максимальные значения среди всех миноров. В каждой подобласти вида $\Omega^{(\nu)}$, одновременно, и выполняется неравенство $0,5 \leq \sqrt{G''_0} \leq 1$, и $\nu$-й минор всюду принимает, по модулю, максимальные значения. Эти подобласти, в общем, могут быть не односвязными. Следующее утверждение следует из следствия к лемме 5:

*Каждая подобласть $\Omega^{(\nu)}$ -замкнутое множество и может быть представлено в виде объединения конечного числа односвязных замкнутых областей, как множество решений в $[0,2^{p-28}]^{24}$ системы полиномиальных неравенств.*

Из сказанных вытекает, что каждая подобласть $\Omega^{(\nu)}$ представляется в виде объединения $\Omega^{(\nu)} = \bigcup_{c \leq T_0} \Omega(\nu,c)$, где $\Omega(\nu,c)$ - односвязные подобласти ($T_0$ не зависит от $p$). Тогда, выделяя, среди всех подобластей ту, которая содержит максимальный кусок $\pi(\nu)$ поверхности $\Pi'_p$, можем писать

$$I_p \leq 2tT_0(2^{28-p})^{4/11} \int_{\pi(\nu)} d\sigma. \quad (13)$$

4) Фиксируем произвольно выбранную подобласть $\Omega(\nu,c)$. В ней система леммы 3 допускает, в произвольной окрестности данного решения, однозначную разрешимость *по одним и тем же переменным*. Эти решения в рассматриваемой подобласти могут быть многолистными. Оценим количество листов разрешимости.

Пусть система (5) допускает однозначную разрешимость, скажем, относительно первых 9-ти переменных $(u_1,v_1,u_2,v_2,u_3,v_3,u_4,v_4,u_5)$, в окрестности некоторого решения $\bar{u}_0$. Оставшиеся переменные, которых мы обозначаем как $\xi_1,...,\xi_{15}$, являются свободными, и пусть $\delta(u_0)$ - область их изменения. Обозначим $\omega(\nu)$ объединение всех областей $\delta(u_0)$, отвечающих всевозможным точкам $\bar{u}_0 \in \pi(\nu)$. Отображение $\varphi: \pi(\nu) \to \omega(\nu)$ такое, что $\varphi(u_0) = \bar{\xi}_0$, $\bar{\xi}_0 = (\xi_1,...,\xi_{15})$ определяет, согласно лемме 1 из [20, с. 538], некоторое $f$-листное накрытие, т. е. система 3 в $\Omega(\nu,c)$ имеет $f$-листную разрешимость. Эта система является полиномиальной, и, поэтому, рассматривая результанты, с помощью последовательного исключения неизвестных из системы, получаем, что $f$ не превосходит некоторого постоянного $T > 0$. Следовательно, область $\Omega(\nu,c)$ можно разбить на не более, чем $T$ подобластей $\Delta_\mu, \mu = 1,..., f, f \leq T$, в каждой из которых система леммы 3 допускает однолистную разрешимость.

5) Обозначим $\pi_1$ такую часть $\pi_1 \subset \pi(\nu)$, которая лежит в области $\Delta_1$ [16, стр. 79]. Система (5), вместе с каждым решением $\bar{u} = (u_1,...,u_{24})$, имеет, также и множество решений вида $\bar{u}' = (u'_1,...,u'_{24})$, с $\bar{u}' = \bar{u} + \bar{a}$ $\bar{a} = (a,b,...,a,b)$, $a,b \in R$. При $\bar{u} = (u_1,...,u_{24}) \in \pi$, точка $\bar{u}'$, с произвольными действительными $a,b \in R$, может лежать на $\pi_1$ если $|a|,|b| \leq 2^{(p-28)/11}$ (иначе, $\bar{u}'$ не принадлежит кубу $[0,(2h)^{-\alpha}]^{24}$). Далее, согласно лемме 4, [16, стр. 79] $G(\bar{u}) = G(\bar{u} + \bar{a})$ при $\bar{u} \in \pi_1$. Множество векторов $\bar{a}$, с $a,b \in R$ образует двумерное подпространство, в $R^{24}$, которого обозначаем, как $V$. В $\Omega(v,c)$, а значит, и в $\Delta_1$, один из миноров матрицы Якоби, принимает, по модулю, максимальные значения. Пусть таким минором будет, например, минор, у которого столбцы получены дифференцированием по первым 9 переменным $u_1, v_1,...$. Тогда, поверхность $\pi_1$ имеет параметрическое представление:

$$u_1 = u_1(\xi_1,...,\xi_{15}),$$
$$\ldots \quad \ldots \quad \ldots$$
$$u_9 = u_9(\xi_1,...,\xi_{15}).$$

Определим на $\pi_1$ отношение эквивалентности, считая $\bar{u} = \bar{u}'$, тогда и только тогда, когда $\bar{u} - \bar{u}' \in V$ (см. [16, стр. 79]). Каждый класс представляет собой линейное подмногообразие вида $\bar{u} + V$, где $\bar{u}$ произвольное решение системы (5). Удобно, также это отношение рассмотреть сначала во всем $R^{15}$. Тогда, все многообразие решений системы (5) можно представить как $\pi_1' + V$, где $\pi_1'$ - решение (его можно отождествить с фактор множеством), проходящее через фиксированную точку $\bar{u}_0$. Если, теперь, мы возьмем подмногообразие $\pi_1' \subset \pi_1$, с фиксированными значениями $\xi_{14} = \xi_{14}^0$, $\xi_{15} = \xi_{15}^0$, размерности 13, с максимальной площадью, то $\pi_1$ покроется при помощи всевозможных параллельных переносов подмногообразия $\pi_1'$ на векторы $\bar{a} \in V$, с $|a|,|b| \leq 2^{(p-28)/11}$. Поэтому, площадь бесконечно малого элемента поверхности $\pi_1$ можно представить в виде $\Delta\bar{a} \cdot \Delta_V \sigma$, где $\Delta\bar{a} = \Delta a \Delta b$, а $\Delta_V \sigma$ означает площадь проекции бесконечно малого элемента поверхности $\pi_1'$ (площадь, которую обозначаем как $\Delta\sigma$) на подпространство, ортогональное к подпространству $V$. Следовательно, $\Delta_V \sigma \leq \Delta\sigma$ и поэтому, согласно (13), имеем

$$I_p \leq 2tT_0(2^{28-p})^{4/11} \int_{\pi(v)} d\sigma \leq 2tT_0(2^{28-p})^{4/11} \times$$

$$\times \int_{|a|,|b| \leq 2^{(p-28)/11}} dadb \int_{\pi_1'} d\sigma \leq 2tT_0(2^{28-p})^{2/11} \int_{\pi_1'} d\sigma. \qquad (14)$$

Заметим, что здесь мы воспользовались теоремой из [18, стр. 221] о повторном интегрировании. Мы внизу будем неоднократно воспользоваться ею замечая, что условия теоремы всегда выполнены, т. к. интегрирование каждый раз будет проводится в замкнутой области, где непрерывная функция равномерно непрерывна.

6) Оценим последний поверхностный интеграл в (14). Из элементов матрицы $A_0 = A_0(\bar{u})$ составим новую блочную матрицу $D_0 = D_0(\bar{u})$ порядка 13 следующим образом:

$$D_0 = \begin{pmatrix} B_1 & B_2 \\ D' & D'' \end{pmatrix},$$

где $B_1$ - максимальный минор порядка 9, $B_2$ - составлена, произвольным образом, из оставшихся столбцов матрицы $A_0$, при этом важно, чтобы все независимые переменные присутствовали в элементах матрицы $D_0$ (поскольку она образована столбцами $A_0$, то все независимые переменные войдут в одночлены наивысших степеней, входящих в $A_0$). Здесь $D'$ -нулевая матрица размера $4 \times 9$, а $D''$ - матрица порядка 4 вида

$$D'' = \begin{pmatrix} 2^{(28-p)/11}\xi_{10} & 0 & 0 & 0 \\ 0 & 2^{(28-p)/11}\xi_{11} & 0 & 0 \\ 0 & 0 & 2^{(28-p)/11}\xi_{12} & 0 \\ 0 & 0 & 0 & 2^{(28-p)/11}\xi_{13} \end{pmatrix}.$$

Имеем:

$$\det D_0 = \det B_1 \cdot 2^{4(28-p)/11}\xi_{10}\xi_{11}\xi_{12}\xi_{13}.$$

Поскольку, миноры матрицы $(B_1 \; B_2)$ являются, также и минорами матрицы $A_0$, то все они по модулю не превосходят 1 и поэтому, из условия

$$\sqrt{G_0''} = \left|\det A_0 \cdot {}^t A_0\right| \leq 1,$$

учитывая границы изменения независимых переменных $\xi_1,...,\xi_{13}$, получаем $|\det D_0| \leq 1$. Следовательно, заменяя условие интегрирования в (12) на $|\det D_0| \leq 1$, мы лишь увеличим поверхностный интеграл. Оценим теперь поверхностный интеграл на правой части (14) при таких условиях. С этой целью, заменим его кратным интегралом по области независимых переменных. Эта область получается проектированием $\pi_1'$ на подпространство независимых переменных $\xi_1,...,\xi_{13}$ (переменные $\xi_{14},\xi_{15}$ фиксированы). Тогда, обозначая область независимых переменных $\tau_0$ и переходя в поверхностном интеграле к независимым переменным, получим следующее неравенство, используя соотношение для элемента площади в лемме 1:

$$\int_{\pi_1'} d\sigma \leq \sqrt{\binom{24}{11}} \int_{\tau_0, |\det D_0| \leq 1} d\xi_1 \cdots d\xi_{13}. \qquad (15)$$

Введем матрицу $D_1$ (также $A_1$), получающуюся расположением элементов столбцов матрицы $D_0$ (соответственно $A_0$), последовательно, в строку, с последующим взятием транспонированной матрицы Якоби полученной системы функций. Эта матрица имеет размер $13\times 169$ (соответственно, размер $24\times 216$). Рассмотрим подобласть $\tau \subset \tau_0$, где $\det(D_1 \cdot {}^t D_1) > G_1$ при некотором $G_1 > 0$ определяемом ниже. В оставшееся части области $\tau_0$ неравенство $|\det D_0| \leq 1$ будет заменяться на условие $\det(D_1 \cdot {}^t D_1) \leq G_1$ и проводиться аналогичная оценка.

7) Далее, будем использовать схему работы [11]. Для объема области $\tau$, сначала заметим, что верна тривиальная оценка $\mu(\tau) \leq 2^{13(p-28)/11}$. Как показано выше, $\det D_0 \neq 0$. Мы оценим более общий интеграл, заменяя условие $|\det D_0| \leq 1$ на $|\det D_0| \leq H$, при этом предполагаем,

$$\int_{\tau, |\det D_0| \leq H} d\xi_1 \cdots d\xi_{13} \leq \sum_{j=1}^{\infty} E_j,$$

где

$$E_j = \int_{2^{-j} H \leq |\det D_0| \leq 2^{1-j} H} d\xi_1 \cdots d\xi_{13}.$$

Пусть $\rho_1 = \rho_1(\overline{\xi'}), \ldots, \rho_{13} = \rho_{13}(\overline{\xi'})$, где $\overline{\xi'} = (\xi_1, \ldots, \xi_{13})$, являются сингулярными числами матрицы $D_0$, $\rho_1 \geq \ldots \geq \rho_{13}$. Тогда, из неравенства

$$2^{-j} H \leq \rho_1 \cdots \rho_{13} \leq \rho_{13} \rho_1^{12}$$

выводим

$$2^{-j} H \rho_1^{-12} \leq \rho_{13}.$$

Полагая $D_0 = (d_{ij})$, по лемме Шура [21, стр. 63], имеем:

$$\rho_1 \leq \max_i \sum_j |d_{ij}| \leq 12 + 72 \cdot 2^{2(p-28)/11}.$$

Следовательно,

$$\rho_{13} \geq 2^{-j} H \left[12 + 72 \cdot 2^{2(p-28)/11}\right]^{-12}. \qquad (16)$$

Перейдем к оценке $E_j, j = 1, 2, \ldots$. Пользуясь интегральным представлением для обратного значения определителя (см. [21, стр. 125, задача 35]) можем писать:

$$\frac{E_j}{2^{1-j} H} \leq \int_{2^{-j} H \leq |\det {}^t D_0| \leq 2^{1-j} H} \frac{d\overline{\xi'}}{|\det {}^t D_0|} \leq c_0 \int_{2^{-j} H \leq |\det {}^t D_0| \leq 2^{1-j} H} d\overline{\xi'} \int_{\|{}^t D_0 \overline{\alpha}\| \leq 1} d\overline{\alpha}, \qquad (17)$$

где $c_0 = \pi^{-13/2} \Gamma(13/2)$.

Далее, из неравенства

$$1 \geq \|D_0 \overline{\alpha}\|^2 = ({}^t D_0 \cdot D_0 \overline{\alpha}, \overline{\alpha}) \geq \rho_{13}^2 \|\overline{\alpha}\|^2 \geq \rho_{13}^2 |\alpha_i|; \|\overline{\alpha}\|^2 = \sum_{i=1}^{13} |\alpha_i|^2$$

для всех $i$, согласно (16), следует оценка

$$|\alpha_i| \leq \left( \sum_{l=1}^{13} |\alpha_l|^2 \right)^{1/2} \leq \rho_{13}^{-1} \leq 2^j H^{-1} \lambda; \quad \lambda = \left[ 12 + 72 \cdot 2^{2(p-28)/11} \right]^{12},$$

для переменных интегрирования во внутреннем интеграле в (17). Введем в рассмотрение шар:

$$K = \{ \overline{\alpha} \mid \|\overline{\alpha}\| \leq 2^j H^{-1} \lambda \}.$$

(17) можно переписать в виде

$$\frac{E_j}{2^{1-j} H} \leq c_0 \int_\tau d\overline{\xi}' \int_{K, \|{}^t D_0 \overline{\alpha}\| \leq 1} d\overline{\alpha} = c_0 \int_{\overline{\alpha} \in K} d\overline{\xi}' \int_{\overline{\xi}' \in \tau, \|{}^t D_0 \overline{\alpha}\| \leq 1} d\overline{\alpha}. \quad (18)$$

Из шара $K$ удалим все полосы $K_m$ ($m = 1, ..., 13$), определяемые условиями

$$|\alpha_m| \leq (1/13) G_1^{-1} 2^{-12j} H^{13} \cdot 2^{13(28-p)/11} \lambda^{-12}, \quad (19)$$

$$|\overline{\alpha}| \leq 2^j H^{-1} \lambda, \; i \neq m. \quad (20)$$

Обозначим $K_0 = \bigcup_{m=1}^{13} K_m$ и оценим меру $K_0$:

$$\mu(K_0) \leq c_0 H G_1^{-1}. \quad (21)$$

На правой части (18), разбивая кратный интеграл на два интеграла, и определяя первый из них условием $\overline{\alpha} \in K_0$, а второй – условием $\overline{\alpha} \in K \setminus K_0$, оценим первый интеграл тривиально, используя найденную оценку:

$$c_0 \int_{\overline{\alpha} \in K_0} d\overline{\alpha} \int_{\overline{\xi}' \in \tau, \|{}^t D_0 \overline{\alpha}\| \leq 1} d\overline{\xi}' \leq c_0 H G_1^{-1}. \quad (22)$$

8) Оценка оставшейся части интеграла (18), где $\overline{\alpha} \in K \setminus K_0$ проводится как в [46] с заменой $s$ на 13. Во внутреннем интеграле на крайней правой части (18), при каждом фиксированном $\overline{\alpha} \in K \setminus K_0$, произведем замену переменных по формулам

$$\overline{\eta} = {}^t D_0(\overline{\xi}') \overline{\alpha}$$

(свободными переменными являются компоненты вектора $\overline{\xi}'$). Матрица Якоби замены переменных равна обратной к матрице

$$J = \frac{\partial(\eta_1, ..., \eta_{13})}{\partial(\xi_1, ..., \xi_{13})} = \left( \frac{\partial {}^t D_0}{\partial \xi_1} \overline{\alpha} \cdots \frac{\partial {}^t D_0}{\partial \xi_{13}} \overline{\alpha} \right),$$

причем $\partial({}^tD_0)/\partial\xi_j$ обозначает матрицу, полученную путем дифференцирования всех элементов матрицы ${}^tD_0(\bar\xi')$ ${}^t$ по $\xi_j$. Имеем:

$$\int\limits_{\bar\alpha\in K\setminus K_0} d\bar\alpha \int\limits_{\bar\xi'\in\tau, \|{}^tD_0\bar\alpha\|\le 1} d\bar\xi' = \int\limits_{\bar\alpha\in K\setminus K_0} d\bar\alpha \int\limits_{\bar\xi'\in\tau, \|\bar\eta\|\le 1} |J|^{-1} d\bar\eta. \quad (23)$$

Левая часть этого равенства конечна. Равенство $|J|=0$ как полиномиальное уравнение имеет множество решений являющимся замкнутым жордановым множеством в $(K\setminus K_0)\times\tau$. При каждом $\bar\xi'\in\tau$ замкнутое жорданово множество, где $|J|=0$ не может содержать точек $\bar\alpha\in\overline{K\setminus K_0}$ вместе с некоторой открытой окрестностью. Действительно, транспонированная матрица ${}^tJ$ имеет вид

$${}^tJ = \begin{pmatrix} {}^t\bar\alpha\cdot\dfrac{\partial{}^tD_0}{\partial\xi_1} \\ \cdots \\ {}^t\bar\alpha\cdot\dfrac{\partial{}^tD_0}{\partial\xi_{13}} \end{pmatrix}.$$

Тогда, для произвольной матрицы $M$ равенство $\bar\alpha\cdot M=0$ может выполняться в открытом подмножестве $\overline{K\setminus K_0}$ только тогда, когда $M=0$. Согласно определению $\tau$, $G_1>0$ и, поэтому, нетривиальная линейная комбинация строк матрицы $D_1$ не равна нулю. Линейная комбинация строк матрицы имеет вид $\bar\alpha\cdot M$, где матрица $M$ равна линейной комбинации матриц $\partial{}^tD_0/\partial\xi_j$ и, следовательно, по сказанному выше, не все ее элементы равны нулю. Поэтому, строки матрицы равенство ${}^tJ$ не могут быть линейно зависимыми в открытом подмножестве $\overline{K\setminus K_0}$. Это означает, что жорданово множество решений уравнения $|J|=0$ не содержит внутренних точек в $(K\setminus K_0)\times\tau$ (действительно, если $V$ такая окрестность, то для любой точки $(\bar\alpha,\bar\xi')\in V, \bar\alpha\in K\setminus K_0, \bar\xi'\in\tau$ найдутся такие окрестности $\bar\alpha\in V', \bar\xi'\in\tau'$, что $V'\times\tau'\subset(K\setminus K_0)\times\tau$), и, поэтому, имеет нулевую жорданову меру. Тогда, мы можем взять некоторое открытое покрытие, состоящее из кубов, с произвольно малым суммарным объемом. Левая часть соотношения (23), взятое по такому покрытию (см. п. 1) и 2)), также, произвольно мала. Если получить равномерную, не зависящую от покрытия, оценку интеграла по замкнутому дополнению этого покрытия, то тем самым мы оценили интеграл, исключением из рассмотрения упомянутое подмножество жордановой меры нуль (т.е. (23) будет справедливым в несобственном смысле).

Итак, оценим интеграл на правой части (23) при таких условиях. При каждом $\overline{\eta}$ обозначим $\tau(\overline{\eta})$ подмножество в $K \setminus K_0$, все точки которого подчинены неравенству $\|{}^t D_0 \overline{\alpha}\| \leq 1$. Рассмотрим матрицу ${}^t J$ как матрицу линейного преобразования, ставящего каждому вектору $\overline{\beta} \in R^{13}$ в соответствие вектор ${}^t J \overline{\beta}$. Очевидно, оно линейно и по $\overline{\alpha}$. Следовательно, имеем билинейное отображение $\Phi : (\overline{\alpha}, \overline{\beta}) \mapsto {}^t J \overline{\beta}$. Для каждой пары $(\overline{\alpha}, \overline{\beta}) \in R^{26}$ выполняется равенство $\Phi(\overline{\alpha}, \overline{\beta}) = D_1(\overline{\alpha} \otimes \overline{\beta})$, при этом, если ${}^t \overline{\alpha} = (\alpha_1, ..., \alpha_{13})$ и ${}^t \overline{\beta} = (\beta_1, ..., \beta_{13})$, то символ ${}^t(\overline{\alpha} \otimes \overline{\beta})$ будет обозначать прямое произведение $(\alpha_1 \beta_1, ..., \alpha_1 \beta_{13}, ..., \alpha_{13} \beta_1, ..., \alpha_{13} \beta_{13})$ (также называемое кронекеровским (см. [22, стр.80] или [23, стр. 235]). Из (23), переставляя порядки интегрирований, получаем:

$$\int\limits_{\|\overline{\eta}\| \leq 1} d\overline{\eta} \int\limits_{\tau(\overline{\eta})} |J|^{-1} d\overline{\alpha} = c_0 \int\limits_{\|\overline{\eta}\| \leq 1} d\overline{\eta} \int\limits_{\tau(\overline{\eta})} d\overline{\alpha} \int\limits_{\|D_1(\overline{\alpha} \otimes \overline{\beta})\| \leq 1} d\overline{\beta}. \qquad (24)$$

Рассмотрим кратный внутренний интеграл по $\overline{\alpha}$ и $\overline{\beta}$:

$$\int\limits_{\tau(\overline{\eta}), \|D_1(\overline{\alpha} \otimes \overline{\beta})\| \leq 1} d\overline{\alpha} d\overline{\beta}. \qquad (25)$$

Пусть сингулярное разложение матрицы $D_1$ имеет вид $D_1 = Q \Sigma T$, где $Q$ и $T$ являются, соответственно, ортогональными матрицами размеров 13 и 169, а $\Sigma$ имеет вид $(\Sigma_1, \Sigma_2)$, $\Sigma_1$ - диагональная матрица, составленная из сингулярных чисел $\sigma_1, ..., \sigma_{13}$ матрицы $D_1$, $\Sigma_2$ - нулевая матрица (конечно столбцы $D_1$ могут располагаться в $\Sigma$ в произвольном порядке, в зависимости от расположения столбцов матрицы $T$).

Рассмотрим интеграл (25) и произведем замену переменных $t_i = \alpha_i \beta_i, i = 1, ..., 13$. Прежде чем применить лемму 1, проведем следующие рассуждения. Из подстановки находим: $\beta_i = t_i \alpha_i^{-1}$. Обозначим это равенство условно в виде $\overline{\beta} = \overline{t} \overline{\alpha}^{-1}$. Тогда можем писать

$$D_1(\overline{\alpha} \otimes \overline{\beta}) = D_1(\overline{\alpha} \otimes \overline{t} \overline{\alpha}^{-1}).$$

Согласно (19-20), для всех точек из $K \setminus K_0$, при всех $i$ выполняются неравенства

$$(1/13) G_1^{-1} 2^{-12j} H^{13} \cdot 2^{13(28-p)/11} \lambda^{-12} \leq \alpha_i \leq 2^j H^{-1} \lambda. \qquad (26)$$

Ниже мы будем накладывать на $G_1(\overline{\xi}')$ условие $G_1(\overline{\xi}') \geq G_1$, где $G_1 > 0$ — та же самая постоянная, рассмотренная выше, не зависящая от $h$ и определяемая точнее позже. Теперь применим лемму 1, производя замену переменных по указанным выше формулам:

$$\int\limits_{\tau(\bar{\eta})} d\bar{\alpha} \int\limits_{\|D_1(\bar{\alpha}\otimes\bar{t}\bar{\alpha}^{-1})\|} d\bar{\beta} =$$

$$= \int d\bar{t} \int\limits_{t_i=\alpha_i\beta_i, \|D_1(t\bar{\beta}^{-1}\otimes\bar{\beta})\|} \frac{ds}{\sqrt{\alpha_1^2+\beta_1^2}\cdots\sqrt{\alpha_{13}^2+\beta_{13}^2}}. \quad (27)$$

Преобразуем поверхностный интеграл в кратный:

$$\int\limits_{t_i=\alpha_i\beta_i, \|D_1(t\bar{\beta}^{-1}\otimes\bar{\beta})\|} \frac{ds}{\sqrt{\alpha_1^2+\beta_1^2}\cdots\sqrt{\alpha_{13}^2+\beta_{13}^2}} \leq \int \frac{d\alpha_1\cdots d\alpha_{13}}{\alpha_1\cdots\alpha_{13}},$$

причем границы изменения $\alpha_i$ определяются неравенствами (26).

Для оценки сверху повторного интеграла на правой части (27) изменим порядки интегрирований:

$$\int \frac{d\alpha_1\cdots d\alpha_{13}}{\alpha_1\cdots\alpha_{13}} \int\limits_{\|D_1(\bar{\alpha}\otimes\bar{t}\bar{\alpha}^{-1})\|\leq 1} dt_1\cdots dt_{13}.$$

Рассмотрим в $R^{169}$ 26-мерное многообразие:

$$t_{11}=\alpha_1\beta_1,\ldots,t_{1,13}=\alpha_1\beta_{13},\ldots,t_{13,1}=\alpha_{13}\beta_1,\ldots,t_{13,13}=\alpha_{13}\beta_{13}.$$

Внутренний интеграл в последнем кратном интеграле можно представить в виде интеграла по поверхности линейного многообразия $\bar{\alpha}\otimes\bar{t}\bar{\alpha}^{-1}$ в $R^{169}$ размерности 13. Элемент площади имеет вид

$$\left|U\cdot{}^tU\right|dt_1\cdots dt_{13},$$

где

$${}^tU = \begin{pmatrix} 1 & 0 & \cdots & 0 & \alpha_2\alpha_1^{-1} & 0 & \cdots & 0 & \cdots \\ 0 & \alpha_1\alpha_2^{-1} & \cdots & 0 & 0 & 1 & \cdots & 0 & \cdots \\ \cdots & \cdots & \cdots & \cdots & \cdots & \cdots & \cdots & \cdots & \cdots \end{pmatrix}.$$

Следовательно,

$$\int\limits_{\|D_1(\bar{\alpha}\otimes\bar{t}\bar{\alpha}^{-1})\|\leq 1} dt_1\cdots dt_{13} = \int\limits_{\|D_1\bar{x}\|\leq 1} \frac{ds}{\left|U\cdot{}^tU\right|},$$

где поверхностный интеграл взят по куску поверхности $\bar{\alpha}\otimes\bar{t}\bar{\alpha}^{-1}$, удовлетворяющей условиям, указанным под знаком интеграла. Матрица $U$ содержит единичную подматрицу, поэтому $\left|U\cdot{}^tU\right|\geq 1$. По неравенству Адамара ([21, стр. 154]) имеем оценку сверху $\left|U\cdot{}^tU\right|\leq 2^{26j+26(p-28)/11}H^{-28}\lambda^{26}G_1^2 = Y$. Преобразуем линейное многообразие $\bar{\alpha}\otimes\bar{t}\bar{\alpha}^{-1}$, действуя на него матрицей $T$ из сингулярного разложения матрицы $D_1$. Так как $T$ ортогональная матрица, то от этого преобразования (т. е. замены переменных $\bar{u}=T\bar{x}$) вид интеграла не изменится. Учитывая вышеприведенные оценки снизу и сверху, получаем:

$$Y^{-1}\int\limits_{\|\Sigma\overline{u}\|\leq 1} d\sigma \leq \int\limits_{\|D_1(\overline{\alpha}\otimes\overline{t}\overline{\alpha}^{-1})\|\leq 1} dt_1\cdots dt_{13} = \int\limits_{\|D_1\overline{x}\|\leq 1} \frac{ds}{|U\cdot {}^tU|} \leq \int\limits_{\|\Sigma\overline{u}\|\leq 1} d\sigma,$$

причем, $d\sigma$ - элемент поверхности на многообразии $\overline{u}=T(\overline{\alpha}\otimes\overline{t}\overline{\alpha}^{-1})$ (эти оценки показывают, что поверхностные интегралы по $ds$ и $d\sigma$ сходятся или расходятся одновременно). Если это многообразие имеет размерность меньшую, чем 13, то вследствие избытка переменных поверхностный интеграл обращается в бесконечность, т. е $|J|=0$. Как сказано выше, равенство $|J|=0$ может выполняться только на множестве тех $\overline{\xi}'$, которые образуют подмножество жордановой меры нуль. Поэтому, можно считать, что рассматриваемое линейное многообразие имеет размерность 13. Обозначим первые компоненты вектора $\overline{u}$ через $u_1,...,u_{13}$. Если подмногообразие, порожденное этими переменными имеет размерность меньшую 13, то опять налицо рассмотренный выше случай. Следовательно, все переменные независимые и поэтому имеем:

$$\int\limits_{\|\Sigma\overline{u}\|} d\sigma = \int\limits_{\sigma_1 u_1^2+\cdots+\sigma_{13}u_{13}^2\leq 1} d\overline{u} = c'\sigma_1^{-1}\cdots\sigma_{13}^{-1} = c'\det(D_1{}^tD_1)^{-1/2} = c'\delta^{-1} \qquad (28)$$

(см. [21, стр. 148]) ($c'$-постоянная). Выполняя обратное преобразование мы находим $T^{-1}\overline{u} = \overline{\alpha}\otimes\overline{t}\overline{\alpha}^{-1}$. Используя границы изменения (26), и интегрируя по $\alpha_j$ под знаком интеграла находим:

$$\int\limits_{13^{-1}G_1^{-1}2^{-12j}H^{13}\lambda^{-13}(2h)^{13\alpha}t_j}^{\lambda H^{-1}2^jt_j} \frac{d\alpha_i}{\alpha_i} = \log(13G_1 2^{13j}H^{-14}\lambda^{13}(2h)^{13\alpha}).$$

Из сказанных выше заключаем:

$$\int\limits_{\|\eta\|\leq 1} d\overline{\eta} \int\limits_{\tau(\overline{\eta})} d\overline{\alpha} \int\limits_{\|D_1(\overline{\alpha}\otimes\overline{\beta})\|\leq 1} d\overline{\beta} \leq c_0\delta_1^{-1}\wp_j^{13}, \qquad (29)$$

где $\wp_j = 1+\log(13G_1 2^{13j}H^{-14}\lambda^{13}(2h)^{-13\alpha})$. Тогда, согласно найденной оценке (29), получаем следующую границу для интеграла по $K\setminus K_0$ на правой части (18):

$$\int\limits_{\overline{\alpha}\in K\setminus K_0} d\overline{\alpha} \int\limits_{\overline{\xi}'\in\tau, \|{}^tD_0\overline{\alpha}\|\leq 1} d\overline{\xi}' \leq c_0^2 H\left(\sum_{j=1}^{\infty}\wp_j^{13}2^{1-j}\right)\delta_1^{-1}.$$

Полагая $X = 13G_1 H^{-14}\lambda^{13}(2h)^{-13\alpha}$, займемся оценкой суммы:

$$\sum_{j=1}^{\infty}\wp_j^{13}2^{1-j} = 2\sum_{j=1}^{\infty}[1+13j\log 2+\log X]^s 2^{-j}.$$

Для оценки этой суммы заметим, что если $1+\log X > 338$, то функция

$$\exp\{13\log(1+13j\log 2+\log X)-0.5j\log 2\}$$

монотонно убывает по $j$. При $1+\log X > 338$ эта функция имеет максимальное значение $\leq 338^{13}$. Поэтому,

$$\sum_{j\geq 0}[1+13j\log 2+\log X]^{13} 2^{-j/2} 2^{-j/2} \leq 4(1+338+\log X)^{13}.$$

Итак,

$$\int_{\tau,|\det D_0|\leq H} d\xi_1\cdots d\xi_{13} \leq 2^{17}\cdot 13^{39} c_0^2 H \delta_1^{-1} \wp^{13}; \wp = 1+338+\log X. \quad (30)$$

Дифференцирование в выражении (28) для матрицы $D_1$ проводится по компонентам вектора $\bar{\xi}'$, которые определяются из системы рассматриваемых уравнений и, поэтому, строки матрицы $D_1(\bar{\xi}')$ являются линейными комбинациями строк матрицы $D_1(\bar{u})$ (см. представление в начале п. 9). Матрица $D_1(\bar{u})$ получается дифференцированием по компонентам $\bar{u}$ и, поэтому, матрицы $D_1(\bar{\xi}')$ и $D_1(\bar{u})$ имеют разные размерности: $D_1(\bar{\xi}')$ имеет размерность $13\times 169$, а $D_1(\bar{u})$ – размерность $24\times 169$. При повторном дифференцировании по $\bar{\xi}'$, возникают сложные выражения, куда входят частные производные зависимых переменных по независимым переменным $\bar{\xi}'$. Поэтому, мы заменим полученные оценки (29-30) оценкой, куда входят подматрицы матрицы $D_1(\bar{u})$, получение которых не осложнены трудностями, указанными выше (т. е. дифференцирование ведется только по независимым переменным вектора $\bar{u}$).

Рассмотрим способ построения матрицы $A_1$. Она получается расположением элементов столбцов матрицы $A_0$ последовательно в строку, и с дальнейшим взятием транспонированной матрицы Якоби полученной системы функций по исходным переменным. Матрица $A_0$ имеет размер $9\times 24$. 96 ее элементов постоянные и потому порождают нулевые столбцы в матрице $A_1$, которая имеет размерность $24\times 216$ (нулевые столбцы исключены). Так как старшая форма 3-й степени данного многочлена содержит обе переменные, то каждая строка матрицы $A_1$ содержит ненулевые элементы, причем все элементы либо постоянные, либо одночлены степени 1:

$$A_1 = \begin{pmatrix} 2 & 0 & 6x_1 & 2y_1 & 0 & 1 & 0 & 2x_1 & 2y_1 & 0 & 0 & \cdots \\ 0 & 1 & 0 & 2x_1 & 2y_1 & 0 & 2 & 0 & 2x_1 & 6y_1 & 0 & \cdots \\ 0 & 0 & 0 & 0 & 0 & 0 & 0 & 0 & 0 & 0 & 2 & \cdots \\ 0 & 0 & 0 & 0 & 0 & 0 & 0 & 0 & 0 & 0 & 0 & \cdots \\ \vdots & \vdots & \vdots & \vdots & \vdots & \vdots & \vdots & \vdots & \vdots & \vdots & \vdots & \ddots \end{pmatrix}.$$

Эта матрица содержит минор равный одночлену с произведением независимых перемен-ных $x_i, y_i, i = 1,...,12$. Следовательно, $G_1$ может обращаться в нуль только на подмножестве жордановой меры нуль. Это дает возможность свести оценку интеграла к оценке подобного несобственного поверхностного интеграла, с множеством особых точек нулевой жордановой меры.

9) Матрицу $A_1(\bar{\xi}')$ можно представить в виде $D(\bar{u}) \cdot A_1(\bar{u})$, где $A_1(\bar{u})$ -матрица, введенная выше, а $D(\bar{u})$-матрица вида:

$$D(\bar{u}) = \begin{pmatrix} 1 & 0 & \cdots & 0 & \varphi_{1,11} & \cdots & \varphi_{1,11} \\ 0 & 1 & \cdots & 0 & \varphi_{2,11} & \cdots & \varphi_{2,11} \\ \vdots & \vdots & \ddots & \vdots & \vdots & \ddots & \vdots \\ 0 & 0 & \cdots & 1 & \varphi_{13,11} & \cdots & \varphi_{13,11} \end{pmatrix} = (E_{13} \mid \Phi),$$

где $E_{13}$ - единичная матрица порядка 13, а $\Phi$ - матрица размера 13×11. Поэтому, произвольный минор, например, минор $M_1$, составленный из первых 13 столбцов матрицы $A_1(\bar{\xi}')$ можно представить в виде блочного определителя

$$\delta_1 = \begin{vmatrix} D(\bar{u}) \cdot A_1^{13}(\bar{u}) & 0 \\ \Psi & E_{11} \end{vmatrix},$$

причем, матрица $A_1^{13}(\bar{u})$ - прямоугольная матрица, составленная из первых 13 столбцов матрицы $A_1(\bar{u})$, $\Phi$- состоит из последних 11 строк матрицы $A_1^{13}(\bar{u})$. Выполняя элементарные преобразования над строками последнего определителя, приведем $\delta_1$ к виду

$$\delta_1 = \begin{vmatrix} \left[A_1^{13}(\bar{u})\right] & -\Phi \\ \Psi & E_{11} \end{vmatrix}, \quad (31)$$

где матрица $\left[A_1^{13}(\bar{u})\right]$ составлена из первых 13 строк матрицы $A_1^{13}(\bar{u})$ так, что две блоки первого столбца определителя $\delta_1$ образуют матрицу $A_1^{13}(\bar{u})$: $A_1^{13}(\bar{u}) = \begin{pmatrix} \left[A_1^{13}(\bar{u})\right] \\ \Psi \end{pmatrix}$. Обозначим матрицу определителя (31) через $D_1$. Разобьем область $\tau$ на две части: в первой выполняется условие $G_1(\bar{\xi}') \geq G_1$, а в оставшейся части области $\tau$ имеем $G_1(\bar{\xi}') \leq G_1$. Обозначая через $\mu_1$ и $\mu_2$ площади соответствующих частей поверхности $\pi_1'$, для (15) получаем:

$$\int_{\pi_1'} d\sigma \leq Z(\mu_1 + \mu_2), \quad (32)$$

где $Z$ не зависит от $p$. Для оценки $\mu_1$ воспользуемся соотношениями (8) и (26-29). Имеем:

$$\mu_1 = 40 \cdot 2^{13} \cdot 13^{40} c_0^2 M T T_0 \left(C_{24}^{11}\right)^{1/2} \left(C_{24 \cdot 11}^{13}\right)^{1/2} H G_1^{-1} \wp^{13}, \quad (33)$$

где $M$ - некоторая постоянная, не зависящая от $p$.

10) Оценку величины $\mu_2$ можно следующим образом свести к оценке подобной (32), уже полученной выше. Сначала выделяем область, где выполнено условие вида $\eta \leq M_1(\bar{\xi}') \leq 2\eta$ (соответствующую площадь обозначаем как $\mu_2^{(1)}$), где $M_1 = M_1(\bar{\xi}')$ минор, содержащий элементы, градиенты которых составляют максимальный минор матрицы $A_2(\bar{\xi}')$. Далее,

$$(2\eta)^{-1}\mu_2^{(1)} \leq \int\limits_{\pi_1',\eta \leq |M_1(\bar{\xi}')| \leq 2\eta} \frac{1}{\delta_1} d\sigma = c_0 \int\limits_{\pi_1',} d\sigma \int\limits_{\|D_1\bar{v}\| \leq 1, \eta \leq |M_1(\bar{\xi}')| \leq 2\eta} dv,$$

при этом интеграл на правой части берется по части произведения $\pi_1' \times R^{24}$, где выполняются наложенные условия на переменные. Поменяв порядки интегрирований, получаем:

$$(2\eta)^{-1}\mu_2^{(1)} \leq 2c_0 \int\limits_{\|D_1\bar{v}\| \leq 1} d\bar{v} \int\limits_{\pi_1'(\bar{v}),\eta \leq |M_1(\bar{\xi}')| \leq 2\eta} dv, \quad (34)$$

причем $\pi_1'(\bar{v})$ - проекция поверхности на $\pi_1'$, отвечающая $\bar{v}$ после перемены порядка интегрирования. Произведем во внутреннем поверхностном интеграле замену переменных по формулам: $\bar{\beta} = D_1\bar{v}$ (заметим, что $D_1$ получается дифференцированием по $\bar{u}$) применяя лемму 2. Получим, в обозначениях этой леммы:

$$\int\limits_{\pi_1',\eta \leq |M_1| \leq 2\eta} \sqrt{G} \frac{d\sigma}{\sqrt{G}} = \int\limits_{\omega(\eta)} |\det Q| \sqrt{G'} \frac{d\sigma}{\sqrt{G'}}, \quad (35)$$

причем $\omega(\eta)$ - прообраз поверхности при указанном отображении и $G' = \det(JQ \cdot {}^tQ^tJ)$; далее, $Q$ - матрица Якоби произведенной замены, которая равна обратной к матрице

$$\frac{\partial \bar{\beta}}{\partial \bar{u}} = \frac{\partial(\beta_1,...,\beta_{24})}{\partial(u_1,...,u_{24})},$$

$J = A_0$, $0.5 \leq \sqrt{G} \leq 1$ (согласно (3)). Как известно, строки матрицы $A_0$ образуют подпространство $M$, ортогональное к подпространству $M'$, натянутому на строки матрицы $D(\bar{u})$. Поэтому, $R^{24} = M \oplus M'$ и, следовательно, элемент объема можно представить в виде $d\bar{w} = d\bar{y}d\bar{z}$. Каждый вектор $\bar{w} \in R^{24}$ единственным образом записывается в виде суммы векторов $\bar{y}$ и $\bar{z}$ из подпространств $M$ и $M'$, при этом $\|Q(\bar{y}+\bar{z})\| \leq \|Q\bar{y}\| + \|Q\bar{z}\|$. Пусть $\bar{y}_1,...,\bar{y}_l$ базис, состоящий из строк матрицы $A_0$, $\bar{z}_1,...,\bar{z}_l$ - базис, состоящий из строк матрицы $D = D(\bar{u})$. Производя замену переменных по формулам $\bar{w} = W\bar{x}$, где

$$W = \begin{pmatrix} A_0 \\ D \end{pmatrix},$$

получаем

$$\left|\det Q^{-1}\right| = c_0 \int\limits_{\|Q\overline{w}\|\leq 1} d\overline{w} = c_0 \int\limits_{\|Q(A_0\overline{u}+D\overline{v})\|\leq 1} \frac{d\overline{u}\,d\overline{v}}{\sqrt{G}\sqrt{\left|D\cdot{}^tD\right|}} \geq$$

$$\geq c_0 \int\limits_{\|QA_0\overline{u}\|\leq 1/2} \frac{d\overline{u}}{\sqrt{G}} \int\limits_{\|QD\overline{v}\|\leq 1/2} \frac{d\overline{v}}{\sqrt{\left|D\cdot{}^tD\right|}}.$$

Следовательно,

$$\left|\det Q^{-1}\right| = c_0 \int\limits_{\|Q\overline{w}\|\leq 1} d\overline{w} \geq$$

$$\geq \frac{\Gamma(13)}{\Gamma(13/2)\Gamma(15/2)} \frac{1}{\sqrt{G}\sqrt{\left|D\cdot{}^tD\right|}\sqrt{G'}\sqrt{\left|{}^tD'QQD\right|}}.$$

Подставляя в (35), находим:

$$\int\limits_{\pi_1'(\overline{v}),\eta\leq|M_1|\leq 2\eta} d\sigma \leq c' \int\limits_{\omega(\eta)} G\sqrt{\left|D\cdot{}^tD\right|}\sqrt{\left|{}^tD'QQD\right|}\,d\sigma'; \quad c' = \frac{\Gamma(13/2)\Gamma(15/2)}{\Gamma(13)}.$$

Матрица $QD$ в каждой точке $\overline{x} \in \Pi_1'$, служащей решением системы, является матрицей Якоби обратного преобразования, т. е. *ее обратная совпадает с матрицей Якоби замены $\overline{\beta} = D_1\overline{v}$, взятой по независимым переменным $\overline{\xi}'$*. Заметим, что элементы матрицы $\Phi$, по модулю, не превосходят 1, и, что $\sqrt{D\cdot D} \leq 11^{13}$, $\sqrt{G} \leq 1$. Тогда, обозначая $C_0$ некоторую постоянную, не зависящую от $p$ и $\eta$, можно вернутся к (18), и далее, при помощи рассуждений п. 8) приходим к соотношению типа (28), уже полученного выше, с заменой $\delta_1$ на $\delta_2$ и $H$ на $\eta$:

$$\mu_2^{(1)} \leq 4c_0\eta \int\limits_{\|D_1\overline{v}\|\leq 1} d\overline{v} \int\limits_{\pi_1'(\overline{v}),\eta\leq|M_1|\leq 2\eta} d\sigma.$$

Этот интеграл оценивается как выше, при этом учитывается соотношения $gG_1^{-1} \ll \left|\det M_1\right| \ll \eta^{-1}$, с положительной $g$, зависящей только от $N$ и $k$. Заменяя $\eta = G_1$ на $\eta/2, \eta/4,...$ а затем, суммируя, мы получаем оценку интеграла (в несобственном смысле, согласно замечанию в конце предыдущего пункта), где выполнено неравенство $\left|\det M_1\right| > 0$:

$$\mu_2 \leq T_1 G_1 G_2^{-1} \wp_1^{13},$$

где $T_1$- положительная постоянная, зависящая только от $N$ и $k$. Заметим, что $\wp_1$- выражение, аналогичное $\wp$. Полагая $G_1 = H^{1/2}G_2^{1/2}$, находим оценку:

$$\int_\delta d\sigma \le ZT_2 \cdot H^{1/2} G_2^{1/2} \wp'^{13},$$

где $T_2$ - зависит только от $N$ и $k$, $\wp' = \max(\wp, \wp_1)$. Заметим, что величины $H, G_1, G_2$ определяются равенствами

$$H = 1, G_1 = G_2^{1/2}.$$

12) Теперь оценим $G_2$ снизу. Матрица $D_2$ составлена из частных производных порядка 3 одночленов $u_1^3, u_1^2 v_1, ...$ по переменным вектора $\bar{u}$. Легко видеть, что в каждой строке содержится хотя бы один ненулевой элемент. При этом, все остальные элементы соответствующего столбца равны нулю. Следовательно, $D_2 \cdot {}^t D_2$ диагональная матрица. Тогда,

$$G_2 = \det\left[D(D_2 \cdot {}^t D_2)^t D\right] \ge 2^{13}.$$

Согласно (26), учитывая очевидное соотношение $\tilde{\wp} \ll \log h^{-1}$, имеем:

$$\int_\delta d\sigma \ll \left(\log h^{-1}\right)^{14}.$$

Из (3), (4) теперь следует

$$\int_\delta d\sigma \ll (2h)^{24\alpha - N\alpha - 1}\left(\log h^{-1}\right)^{14}.$$

Так как $h = \sqrt{G} \cdot 2^{-p} \le 72^4 \cdot 2^{-p}$, то ряд

$$\sum_{p=1}^\infty I_p \ll \sum_{p=1}^\infty 2^{-p(18\alpha - 1)} p^{13}$$

а вместе с ним и особый интеграл сходятся. Доказательство теоремы завершено.

### Литература